\documentclass[letterpaper, 10 pt, conference]{ieeeconf}  % Comment this line out if you need a4paper

\IEEEoverridecommandlockouts                              % This command is only needed if 
\usepackage{amsmath,amsfonts,cite,times,url,hyperref,bm,graphicx}
\newtheorem{lemma}{Lemma}

\title{\LARGE \bf
Auxiliary Signal-Based Distance Protection in \\ Inverter-Dominated Power Systems
}

\author{Josh A. Taylor$^{1}$ and Alejandro D. Dom\'inguez-Garc\'ia$^{2}$% <-this % stops a space
\thanks{$^{1}$Josh A. Taylor is with the Department of Electrical and Computer Engineering,
       New Jersey Institute of Technology, Newark, NJ, USA. E-mail:
        {\tt\small jat94@njit.edu}}%
\thanks{$^{2}$Alejandro D. Dom\'inguez-Garc\'ia is with the Department of Electrical and Computer Engineering, University of Illinois at 
	Urbana-Champaign, Urbana, IL, USA. E-mail:
        {\tt\small aledan@illinois.edu}}%
}

\begin{document}

\maketitle
\thispagestyle{empty}
\pagestyle{empty}

%%%%%%%%%%%%%%%%%%%%%%%%%%%%%%%%%%%%%%%%%%%%%%%%%%%%%%%%%%%%%%%%%%%%%%%%%%%%%%%%
\begin{abstract}

Power system protection schemes today rely on currents rising by several orders of magnitude when faults occur. In inverter-dominated power systems, a fault current might be just a few percent larger than normal, making fault detection more difficult. One solution is for the inverter to slightly perturb its output current and/or voltage, i.e., to inject an auxiliary signal, so as to make the system's behavior under faults easier to distinguish from normal. In this paper, we optimize auxiliary signals for fault detection with distance relays. We begin with a standard auxiliary signal design problem for generic static systems. We use duality to reformulate the problem as a bilinear program, which we solve using the convex-concave procedure. We implement the framework in an example based on distance protection, in which the auxiliary signal is negative sequence current.

\end{abstract}

%%%%%%%%%%%%%%%%%%%%%%%%%%%%%%%%%%%%%%%%%%%%%%%%%%%%%%%%%%%%%%%%%%%%%%%%%%%%%%%%
\section{INTRODUCTION}\label{sec:intro}

In power networks supplied by synchronous generators, large currents provide clear information about the existence and location of faults~\cite{Anerson1999psp,phadke2009computer}. In systems supplied primarily by inverter-interfaced resources (IBRs), inner-loop current controllers prevent large currents even during faults, so that a fault current might be only a few percent larger than normal~\cite{hooshyar2017microgrid}. As a result, fault detection schemes for synchronous machine-fed faults can fail in inverter-dominated grids~\cite{banaiemoqadam2019control}.

One way to make faults easier to detect in inverter-fed grids is to slightly perturb the output current when there is suspicion of a fault~\cite{banaiemoqadam2019control,banaiemoqadam2020comprehensive}. This can be accomplished, for example,  by adding new inputs to the inverter's current and/or voltage controllers. For example, in~\cite{saleh2020protection}, when the inverter's terminal voltage is low, the inverter injects harmonic currents. Similarly, the  IEEE Standard 2800 now stipulates that inverter-interfaced resources must inject negative sequence current to make them behave more like synchronous generators during unbalanced faults~\cite{IEEE28002022}. The downside of such schemes is that perturbing the current lowers power quality, though only when the presence of a fault is suspected. The benefit is that no additional hardware is necessary, which can make them compatible with some existing protection setups~\cite{banaiemoqadam2019control}. 

In this paper, we pose the design of such perturbations as an auxiliary signal problem. An auxiliary signal is a perturbation to the system's inputs that facilitates fault detection~\cite{campbell2015auxiliary}. Because we are focusing on phasor-based protection schemes, here we  only consider static system settings. We use duality to reformulate the standard static auxiliary signal design problem as a bilinear program. An advantage of this approach is that we can incorporate general inequality constraints, which, for instance, might encode an inverter's voltage or current limits. Though nonconvex, the problem is not large in scale, and could be solved in a variety of different ways. Here we use the convex-concave procedure (CCP)~\cite{yuille2003concave,lipp2016variations} because it enables the use of convex solvers and is relatively easy to implement.

We now review  the relevant literature. In several recent papers, the inverter injects a current perturbation to make faults easier to detect. As mentioned above,~\cite{saleh2020protection} injects harmonics, and~\cite{banaiemoqadam2019control,banaiemoqadam2020comprehensive} modulate the positive and negative sequence currents. We also note that negative sequence current has been used to assist islanding detection~\cite{karimi2008negative}. There is an established literature stream on auxiliary signal design, of which the textbook~\cite{campbell2015auxiliary} provides comprehensive coverage. The most relevant paper, which is slightly more recent and from which we take direct inspiration, is~\cite{blanchini2017active}, wherein fault detection in linear time-invariant systems is reformulated using Hahn-Banach duality. The present work differs in that we focus on static systems; the resulting bilinear program and its numerical solution are new. Aside from the first author's recent work in~\cite{pirani2022fault}, none of the formalisms in the auxiliary signal literature have been applied to fault detection in power systems.

The original contributions of this papers are as follows.
\begin{itemize}
\item[\bf C1.] We use duality to reformulate fault detection in static systems as a bilinear program. This allows for inequality constraints in the system model, and solution via a variety of tools from nonlinear programming.
\item[\bf C2.]  We use the  CCP to obtain locally optimal solutions of the bilinear programs.
\item[\bf C3.]  We apply the procedure to distance protection. Specifically, we optimize negative sequence current injection to assist detection of a phase-to-phase fault.
\end{itemize}

The remainder of the paper is organized as follows. In Section~\ref{sec: formulation}, we describe the auxiliary signal design problem for static systems. In Section~\ref{sec:duality}, we use duality to reformulate fault detection in static systems as a bilinear program. In Section~\ref{sec:CCP}, we describe its numerical solution via CCP. In Section~\ref{sec:distance}, we apply the procedure to distance protection. We describe several directions of future work in Section~\ref{sec:conclusion}.

\section{PROBLEM FORMULATION} \label{sec: formulation}
We adapt the formulation in Section 3.2 of~\cite{campbell2015auxiliary}. Consider a linear static system that can operate in one of two modes: normal and faulty, which we index by $k=0$ and $k=1$, respectively. Let $x$ denote a vector of system variables that are known a priori or can be measured. Let $\theta$ denote the auxiliary signal: a vector of system variables whose value can be manipulated so as to determine in which mode the system is operating. Then, for $k=0,1$, the system behavior can be described by
\begin{subequations}
\label{sysnucon}
\begin{align}
\Theta_k\theta + X_kx &= H_k \lambda_k,\label{sys1}\\
A_kx &\leq b_k,\label{sys2}
\end{align}
where $\lambda_k$ is unknown but must satisfy
\begin{align}
\left\|\lambda_k\right\|< 1,\quad k=0,1,\label{nucon}
\end{align}
\end{subequations}
 $\Theta_k$, $X_k$,  $H_k$ are matrices of appropriate dimensions, and $b_k$ is a vector of appropriate dimension.

The  relation in \eqref{sys1} can be a result of physical laws, e.g., Kirchhoff's laws. The inequality constraint in \eqref{sys2} describes specific limits that the vector of known/measured variables must satisfy, e.g., the inverter's maximum output current. The unknown vector, $\lambda_k$, describes measurement errors or unmodeled behavior. For example, the actual relation between the auxiliary signal, $\theta$, and the vector of known/measured variables, $x$, might be nonlinear. In this case, $\lambda_k$ could capture higher-order terms in $x$ and $\theta$ that were omitted from~\eqref{sys1}.

Given \eqref{sys1} -- \eqref{sys2}, we ask if there is a reasonable level of noise, $\lambda_0$ or $\lambda_1$, under which the known quantities could have come from model $k=0$ or $k=1$. For intuition, suppose first there is no auxiliary signal ($\theta=0$). If there is a $\lambda_0$ that satisfies (\ref{sysnucon}) for $k=0$ and no feasible $\lambda_1$ for $k=1$, then we can conclude the system is operating normally. If vice versa, then the system is operating in its faulty mode. If there are $\lambda_0$ and $\lambda_1$ such that (\ref{sysnucon}) is feasible for $k=0$ and $k=1$, then we cannot determine if the system has failed or not. In this case, we design~$\theta$ so that (\ref{sysnucon}) is only feasible for $k=0$ or $k=1$, but not both. Next, we derive a new way of doing so.

\section{DUALITY-BASED AUXILIARY SIGNAL DESIGN} \label{sec:duality}
% \todo{below define $\lambda$ or write $\lambda_0$ and $\lambda_1$ belon min}Consider the optimization
Consider the optimization
\begin{subequations}
\label{sigmav}
\begin{align}
\mathcal{P}_0:\quad\min_{x,\omega,\lambda} \quad& \omega\\
\textrm{subject to} \quad& \Theta_k\theta + X_kx = H_k\lambda_k\label{sigmav1}\\
&A_kx \leq b_k\label{sigmav3}\\
&\left\|\lambda_k\right\|^2 \leq \omega,\quad k=0,1,\label{sigmav2}
\end{align}
\end{subequations}
where $\lambda=\left[\lambda_0^{\top},\lambda_1^{\top}\right]^{\top}$. Let $\sigma(\theta)=\omega^*$ denote the optimal objective. Observe that (\ref{sigmav2}) will bind for $k=0$ or $k=1$.
% \todo{Explain that one of the constraints in (2d) will be binding and explain why (even if it might be obvious to some) }
If $\omega^*\geq1$, then at least one of the models is infeasible. If $\left\|\lambda_0\right\|<1$ and $\left\|\lambda_1\right\|\geq1$, or vice versa, then we can  determine if the system is operating normally ($k=0$) or has a fault ($k=1$). We want to design the auxiliary signal, $\theta$, so that this is the case.

We want $\theta$ to be minimal in some sense, which we quantify by the cost $\theta^{\top}Q\theta$, where $Q\succeq0$. We seek to solve
\begin{align}
\mathcal{P}_1:\quad \min_{\theta} \quad \theta^{\top}Q\theta\quad \textrm{subject to} \quad 1\leq\sigma(\theta).\label{optaux1}
\end{align}
In other words, we want to find the smallest auxiliary signal that makes normal and faulty operation distinguishable. We expect the constraint to bind at the optimal solution, so that only one of the models is infeasible.

The constraint is difficult to handle because it contains a minimization. We can eliminate the minimization by replacing $\mathcal{P}_0$ with its dual, a maximization. Because $\mathcal{P}_0$ is a convex quadratic program, strong duality ensures that it has the same objective as its dual (assuming some constraint qualification holds). We may simply constrain the objective to be greater than one, which ensures the maximum objective is at least one.

Replacing $\mathcal{P}_0$ with its dual and dropping the maximum, $\mathcal{P}_1$ becomes
\begin{subequations}
\label{optaux2}
\begin{align}
\mathcal{P}_2:\min_{\theta,\alpha,\beta,\epsilon,\delta} \quad& \theta^{\top}Q\theta\\
\textrm{subject to} \quad& 1\leq\sum_{k=0}^1\beta_k^{\top}\Theta_k\theta -\epsilon_k^{\top}b_k - \delta_k\label{bilinear}\\
& \alpha_0+\alpha_1=1\label{P2c}\\
&\sum_{k=0}^1\beta_k^{\top}X_k + \epsilon_k^{\top}A_k=\bm{0}\\
& 4\alpha_k\delta_k\geq \left\|H_k^{\top}\beta_k\right\|^2,\quad k=0,1\label{SOC}\\
&\alpha_k\geq0,\;\epsilon_k\geq0,\quad k=0,1,\label{P2f}
\end{align}
\end{subequations}
where $\alpha=\left[\alpha_0^{\top},\alpha_1^{\top}\right]^{\top}$, and $\beta$, $\delta$, and $\epsilon$ are similarly defined. Here $\delta$ is a dummy variable that simplifies (\ref{bilinear}), $\beta$ is the dual variable of (\ref{sigmav1}), $\epsilon$ of (\ref{sigmav3}), and $\alpha$ is the dual variable of (\ref{sigmav2}). The following lemma summarizes the result.

\begin{lemma}
If strong duality holds for $\mathcal{P}_0$, then $\mathcal{P}_1$ is equivalent to $\mathcal{P}_2$.
\end{lemma}

By solving $\mathcal{P}_2$, we obtain an auxiliary signal, $\theta$, for which $\theta^{\top}Q\theta$ is minimal and (\ref{sysnucon}) infeasible. If we assume that one of the models is true, $k=0$ or $k=1$, and that the noise is bounded, this guarantees that we can distinguish between normal and faulty operation.

\section{SOLUTION VIA THE CONVEX-CONCAVE PROCEDURE}\label{sec:CCP}
There are several options for solving $\mathcal{P}_2$. We use the CCP because it allows us to use convex solvers and is straightforward to implement~\cite{yuille2003concave,lipp2016variations}. We now describe its implementation for $\mathcal{P}_2$.

$\mathcal{P}_2$ is nonconvex due to the bilinear constraint (\ref{bilinear}). Constraint (\ref{SOC}) is a hyperbolic constraint, which has a second-order cone representation, and the others are linear.

Define
\[
\Psi_k=\begin{bmatrix}
\bm{0} & \Theta_k\\
\Theta_k^{\top} & \bm{0}
\end{bmatrix},
\]
where $\bm{0}$ and $\bm{I}$ respectively denote zero and identity matrices of appropriate dimension. Let $\psi_k$ denote the largest eigenvalue of $\Psi_k$. Then, each
bilinear term in (\ref{bilinear}) can be written
\[
\frac{1}{4}\begin{bmatrix}
\beta_k\\\theta
\end{bmatrix}^{\top}\left(\Psi_k-\psi_k\bm{I}\right)\begin{bmatrix}
\beta_k\\\theta
\end{bmatrix}
+
\frac{1}{4}\begin{bmatrix}
\beta_k\\\theta
\end{bmatrix}^{\top}\left(\Psi_k+\psi_k\bm{I}\right)\begin{bmatrix}
\beta_k\\\theta
\end{bmatrix}.
\]
The first term is concave and the latter convex. We linearize the latter because it is on the right side of the inequality. The linearization is
\begin{align*}
J_k^z(\beta_k,\theta)=&\frac{1}{4}\begin{bmatrix}
\beta_k^z\\\theta^z
\end{bmatrix}^{\top}\left(\Psi_k+\psi_k\bm{I}\right)\begin{bmatrix}
\beta_k^z\\\theta^z
\end{bmatrix}+\\
& \frac{1}{2}\begin{bmatrix}
\beta_k^z\\\theta^z
\end{bmatrix}^{\top}\left(\Psi_k+\psi_k\bm{I}\right)\left(\begin{bmatrix}
\beta_k\\\theta
\end{bmatrix}-\begin{bmatrix}
\beta_k^z\\\theta^z
\end{bmatrix}\right),
\end{align*}
where $\beta_k^z$ and $\theta^z$ are the values of $\beta_k$ and $\theta$ at the $z^{\textrm{th}}$ iteration of the CCP.

The CCP for $\mathcal{P}_2$ is as follows. Choose algorithm parameters $\gamma^0>0$, $\gamma_{\max}>\gamma^0$, $\zeta>1$, $\theta^0$, and $\beta^0$ and set the iteration counter to $z=0$. Then repeat the below steps until a stopping criterion, e.g., convergence of the objective or variables, is satisfied.
\begin{enumerate}
\item Solve the optimization
\begin{align*}
\left(\theta^{z+1},\beta^{z+1}\right)&=\underset{\theta,\alpha,\beta,\epsilon,\delta,\xi}{\textrm{argmin}} \quad \theta^{\top}Q\theta+\gamma^z\xi^{\top}\bm{1}\\
\textrm{subject to}\quad&\xi\geq0\\
&1-\xi\leq \sum_{k=0}^1J_k^z(\beta_k,\theta)-\epsilon_k^{\top}b_k -\delta_k\\
&(\ref{P2c})-(\ref{P2f}).
\end{align*}
\item Set $\gamma^{z+1}=\min\{\zeta\gamma^z,\gamma_{\max}\}$.
\item Set $z=z+1$.
\end{enumerate}
Here $\xi$ is a slack variable that allows $\theta^0$ and $\beta^0$ to be infeasible. The additional term in the objective is a penalty on $\xi$. Each iteration of the CCP is a second-order cone program, which can be solved efficiently using existing software.

\section{APPLICATION TO DISTANCE PROTECTION}\label{sec:distance}
A distance relay measures local current, $i$, and voltages, $e$ (from each phase to ground or between phases). Let $\bar{i}$ and  $\bar{e}$ denote complex numbers, referred to as phasors, associated with $i$ and $e$, respectively. The relay computes an impedance  by dividing the voltage phasor by the current phasor: $\bar{z}=\bar{e}/\bar{i}$. The presence of a fault affects the voltage and current, and therefore the impedance computed by the relay. If the impedance is in the relay's zone of operation, e.g., a circle or quadrilateral on the complex plane, it initiates protective action, i.e., opening a circuit breaker. The zone of operation is based on estimates of the impedance the relay will see under normal conditions and during a fault~\cite{horowitz2022power}. A distance relay can misdiagnose an inverter-fed fault because the current it measures is little different from normal. We remark that this description of distance relaying is simplistic~\cite{roberts1993z}, but serves as context for the following example. 

\begin{figure}[h]
		\centering
\includegraphics[width=0.4\textwidth]{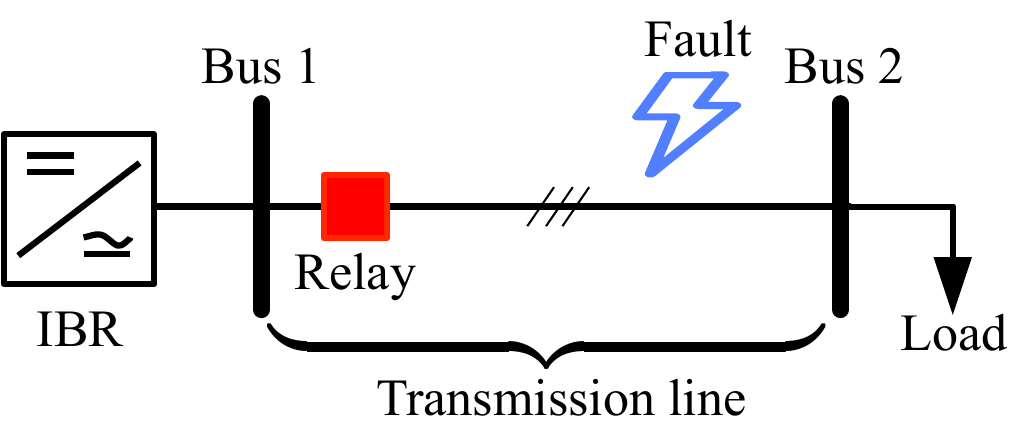}
		\vspace{-0.14in}
	\caption{Single-IBR single-load system. A fault occurs somewhere along the length of the transmission line.}
		\label{fig:fault}
	\end{figure}

\subsection{Single-IBR single-load system}\label{sec:sisl}
Consider an IBR supplying power to a load through a transmission line, as shown in Figure~\ref{fig:fault}.
%; we refer to this system as a single-IBR single-load  system. 
A phase $a$-to-phase $b$ fault occurs somewhere on the line. A distance relay on the line must use local voltage and current measurements to determine if a fault has occurred.

% \begin{figure}[!h]
% \centering
% \begin{tikzpicture}
%     \node [ball] (S){Inverter};
%     \node [block, right=1.5cm of S] (R){Relay};
%     \node [fault, right=.5cm of R] (F) {Fault};
%     \node [ball, right=.5cm of F] (L) {Load};
%     \draw [-, thick] (S) -- node[above]{Line}(R);
%     \draw [-, thick] (R) -- node[above]{}(F);
%     \draw [-, thick] (F) -- node[above]{}(L);
% \end{tikzpicture}
% \caption{An inverter-interfaced source sends power through a line to a load. A fault occurs somewhere along the length of the line.}
% \label{fig:fault}
% \end{figure}

Suppose the relay measures the positive and negative sequence voltage and current phasors, which we denote by $\bar{e}_+$, $\bar{e}_-$, $\bar{i}_+$,  and $\bar{i}_-$, respectively. When there is no fault, the measurements will satisfy
\[
\bar{e}_-\approx \bar{z}_-\bar{i}_-\approx0,\quad \bar{e}_+\approx \bar{z}_+\bar{i}_+,
\]
where $\bar{z}_+$ and $\bar{z}_-$ are estimates of the apparent positive and negative sequence impedances of the line and load. Assume now there is a fault between phases $a$ and $b$. Then measurements will satisfy
\[
\bar{e}_-\approx \bar{z}_{\textrm{f}}\bar{i}_-,\quad \bar{e}_+\approx \bar{z}_{\textrm{f}}\bar{i}_+,
\]
where $\bar{z}_{\textrm{f}}$ is an estimate of the impedance under a fault (see, e.g., Section 5.4.1 of~\cite{horowitz2022power}). These relations are approximate due to uncertainty in the impedances and noisy measurements, which we represent explicitly below.

In this scenario, a synchronous machine would act like a voltage source, injecting large currents with nonzero negative sequence. On the other hand, an IBR behaves like a current source, preventing large currents and keeping $\bar{i}_-\approx0$ even under a fault, leading to an inaccurate estimate of the fault impedance. Similar to~\cite{banaiemoqadam2019control}, we remedy this by having the inverter inject negative sequence current.

Here the auxiliary signal, $\bar{\theta}$, is a phasor. Because it is negative sequence current, it takes the place of $\bar{i}_-$ (which is roughly zero). When there is no fault, the system model is
\begin{subequations}
\label{faultseqa}
\begin{align}
\bar{e}_-&= \bar{z}_{-}\bar{\theta}+ \bar{\lambda}_{-,0}\\
\bar{e}_+&= \bar{z}_+\bar{i}_++\bar{\lambda}_{+,0}.
\end{align}
\end{subequations}
This corresponds to (\ref{sys1}) when $k=0$. When there is a fault, the system model is
\begin{subequations}
\label{faultseqb}
\begin{align}
\bar{e}_-&= \bar{z}_{\textrm{f}}\bar{\theta}+\bar{\lambda}_{-,1}\label{faultseqc}\\
 \bar{e}_+&= \bar{z}_{\textrm{f}}\bar{i}_++\bar{\lambda}_{+,1}.\label{faultseqd}
\end{align}
\end{subequations}
This corresponds to (\ref{sys1}) when $k=1$ (faulty operation). The noises must satisfy
\begin{align}
\left\|\begin{bmatrix}\bar{\lambda}_{-,0}\\\bar{\lambda}_{+,0}\end{bmatrix}\right\|\leq 1,\quad\left\|\begin{bmatrix}\bar{\lambda}_{-,1}\\\bar{\lambda}_{+,1}\end{bmatrix}\right\|\leq 1.\label{faultseqe}
\end{align}
We want to design $\bar{\theta}$ so that either (\ref{faultseqa}) or (\ref{faultseqb}) is feasible, but not both.

In Appendix~\ref{app:real}, we convert (\ref{faultseqa})-(\ref{faultseqe}) to a real-valued system and put it in the form of (\ref{sysnucon}). We can thus fully parameterize the auxiliary signal design optimization, $\mathcal{P}_2$. We remark that it would be more natural here to add the uncertainty directly to the impedances. This would result in multiplicative uncertainty, a topic of future work.

\subsection{Numerical results}
To solve $\mathcal{P}_2$ numerically with the CCP, we must specify the values of $\bar{z}_-,\bar{z}_+,$ and $\bar{z}_{\textrm{f}}$. We set
\begin{align*}
\bar{z}_-&=30 + j35,\\
\bar{z}_+&=30 + j35,\\
\bar{z}_{\textrm{f}}&=26+jx_{\textrm{f}},
\end{align*}
where $x_{\textrm{f}}$, the reactance of the fault path, is a parameter we vary.

All optimizations were carried out in Python using CVXPy~\cite{diamond2016cvxpy} and the solver Gurobi~\cite{gurobi}. All figures were made with Matplotlib~\cite{hunter2007matplotlib}. Each optimization took under a second to solve.

\begin{figure}[t]
\centering
\includegraphics[width = \columnwidth]{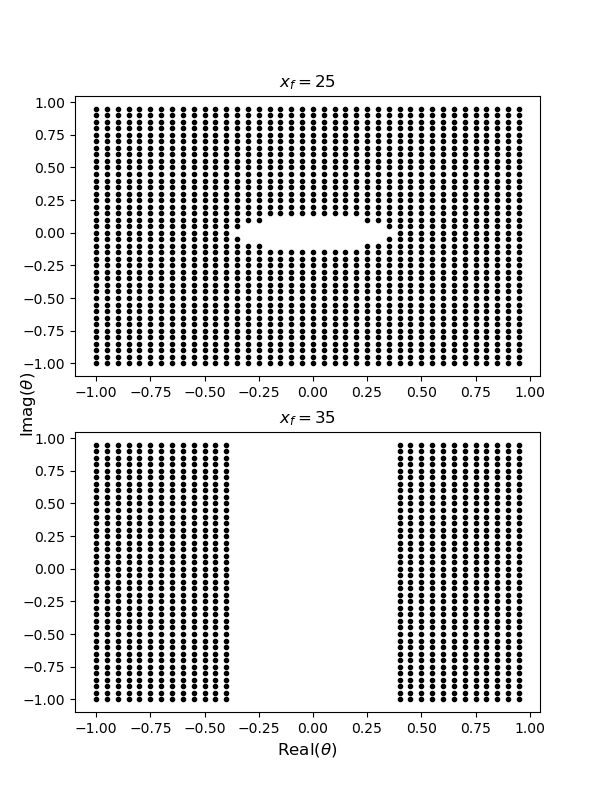}
%\vspace{-0.4in}
\caption{Feasible values of $\bar{\theta}$ (dotted) on the complex plane for $x_{\textrm{f}}=25$ (top) and $x_{\textrm{f}}=35$ (bottom).}
\label{fig:ABfeas}
\vspace{0.2in}
\end{figure}

We first explore the feasibility of $\mathcal{P}_2$. Testing the feasibility of a particular auxiliary signal entails solving $\mathcal{P}_2$ for a fixed value of $\bar{\theta}$. Fixing $\bar{\theta}$ makes the bilinearity linear, which makes $\mathcal{P}_2$ a second-order cone program. Figure~\ref{fig:ABfeas} shows feasible values of $\bar{\theta}$ on the complex plane for two cases: $x_{\textrm{f}}=25$ and $x_{\textrm{f}}=35$. The latter case is more difficult because the fault impedance is closer to the normal impedance seen by the relay. As a result, the set of infeasible auxiliary signals, i.e., which do not make the fault distinguishable form normal operation, is larger when $x_{\textrm{f}}=35$.

\begin{figure}[t]
\centering
\includegraphics[width = \columnwidth]{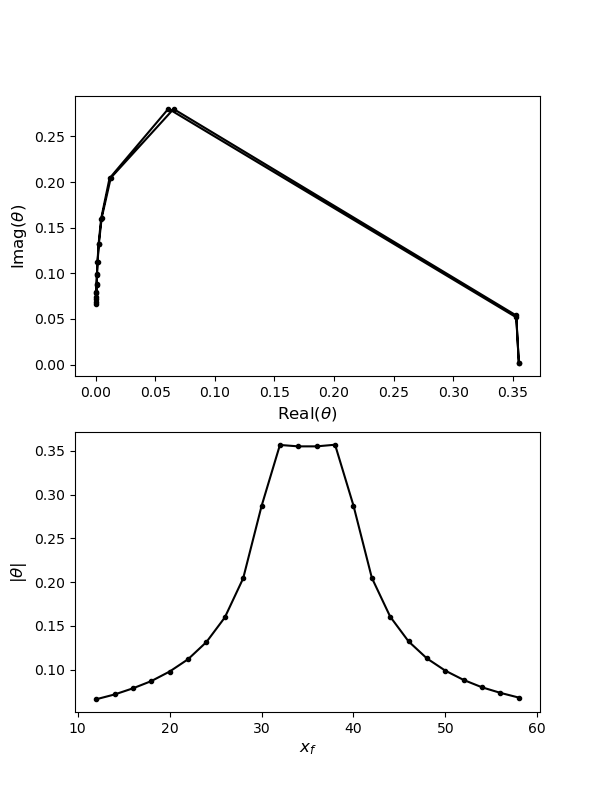}
%\vspace{-0.4in}
\caption{$\bar{\theta}$ (top) and $\left|\bar{\theta}\right|$ (bottom) for $x_{\textrm{f}}\in[12,58]$.}
\label{fig:ABtheta}
\vspace{0.2in}
\end{figure}

We now use the CCP to solve $\mathcal{P}_2$ for the optimal auxiliary signal, $\bar{\theta}$. Figure~\ref{fig:ABtheta} shows $\bar{\theta}$ as $x_{\textrm{f}}$ varies from 12 to 58. The top plot shows the real and imaginary components of $\bar{\theta}$ for each value of $x_{\textrm{f}}$. The bottom shows $\left|\bar{\theta}\right|$ as a function of $x_{\textrm{f}}$. When $x_{\textrm{f}}$ is either very small or very large, it is significantly different from the apparent reactance during normal conditions, making faults easier to detect. As a result, $\bar{\theta}$ is on the bottom left of the top plot, and has small magnitude on the bottom plot. When $x_{\textrm{f}}\approx35$, the fault is more difficult to detect, and a larger auxiliary signal is necessary. In this case $\bar{\theta}$ is on the bottom right of the top plot, and has larger magnitude on the bottom.

The optimization $\mathcal{P}_2$ is a generic framework for obtaining auxiliary signals for static systems. This example demonstrates that it can be solved efficiently, and produces intuitive results when applied to a basic instance of fault detection.

\section{CONCLUSION}\label{sec:conclusion}
We have used duality to reformulate auxiliary signal design for static systems as a bilinear program. The benefit is more modeling flexibility and amenability to the large set of solution techniques for nonlinear programming. Here we have used the CCP because it is easy to implement and relies on convex solvers.

The motivation for this work is fault detection in inverter-fed grids. The problem is difficult because the fault currents of IBRs are smaller than those of synchronous machines. A relatively recent solution is to add perturbations that make faults easier to distinguish from normal conditions. We have formalized this strategy as an auxiliary signal design problem, and applied the resulting procedure in an example with a phase-to-phase fault.

We believe this is a natural approach to fault detection in inverter-dominated power systems, which can be applied in a range of scenarios. Two immediate directions of future work are generalizing the procedure to more than two models, and applying it to the full set of faults seen in three-phase systems. We also intend to apply auxiliary signal formalisms for dynamic systems, e.g.,~\cite{blanchini2017active}, to time-domain protection.

\section*{APPENDIX}\label{app:real}

We put the system in Section~\ref{sec:distance}, (\ref{faultseqa})-(\ref{faultseqe}), in the form of (\ref{sysnucon}), which enables us to parameterize the auxiliary signal design optimization, $\mathcal{P}_2$.

To this end, we convert all quantities from complex- to real-valued. Let
\[
\bar{\theta}= \mu + j \nu
\]
and
\begin{align*}
\bar{e}_-&= f_- + jg_-\\
\bar{i}_- &= c_- + jd_-\\
\bar{z}_-&= r_- + jx_-.
\end{align*}
We replace the subscript with $+$ and $\textrm{f}$, respectively, to denote positive sequence and fault values. Let
\[
M_-=\begin{bmatrix}r_-&-x_-\\r_-&x_-\end{bmatrix},
\]
and define $M_+$ and $M_{\textrm{f}}$ similarly. Let $\mathcal{R}\left[\bar{y}\right]$ and $\mathcal{I}\left[\bar{y} \right]$ denote the real an imaginary parts of $\bar{y}$. Then (\ref{faultseqa}) takes the form
\begin{align*}
\begin{bmatrix}f_-\\g_-\end{bmatrix}&=M_-\begin{bmatrix}\mu\\\nu\end{bmatrix}+\begin{bmatrix}\mathcal{R}\left[\bar{\lambda}_{-,0}\right]\\\mathcal{I}]\left[\bar{\lambda}_{-,0}\right]\end{bmatrix}\\
\begin{bmatrix}f_+\\g_+\end{bmatrix}&=M_+\begin{bmatrix}c_+\\d_+\end{bmatrix}+\begin{bmatrix}\mathcal{R}\left[\bar{\lambda}_{+,0}\right]\\\mathcal{I}\left[\bar{\lambda}_{+,0}\right]\end{bmatrix},
\end{align*}
(\ref{faultseqb}) takes the form
\begin{align*}
\begin{bmatrix}f_-\\g_-\end{bmatrix}&=M_f\begin{bmatrix}\mu\\\nu\end{bmatrix}+\begin{bmatrix}\mathcal{R}\left[\bar{\lambda}_{-,1}\right]\\\mathcal{I}\left[\bar{\lambda}_{-,1}\right]\end{bmatrix}\\
\begin{bmatrix}f_+\\g_+\end{bmatrix}&=M_f\begin{bmatrix}c_+\\d_+\end{bmatrix}+\begin{bmatrix}\mathcal{R}\left[\bar{\lambda}_{+,1}\right]\\\mathcal{I}\left[\bar{\lambda}_{+,1}\right]\end{bmatrix},
\end{align*}
and (\ref{faultseqe}) takes the form
\begin{align*}
&\mathcal{R}\left[\bar{\lambda}_{-,0}\right]^2+\mathcal{I}\left[\bar{\lambda}_{-,0}\right]^2+\mathcal{R}\left[\bar{\lambda}_{+,0}\right]^2+\mathcal{I}\left[\bar{\lambda}_{+,0}\right]^2\leq 1\\
&\mathcal{R}\left[\bar{\lambda}_{-,1}\right]^2+\mathcal{I}\left[\bar{\lambda}_{-,1}\right]^2+\mathcal{R}\left[\bar{\lambda}_{+,1}\right]^2+\mathcal{I}\left[\bar{\lambda}_{+,1}\right]^2\leq 1.
\end{align*}

We can now parameterize (\ref{sysnucon}). With a slight abuse of notation, arrange the variables as
\[
\theta = \begin{bmatrix}\mu\\\nu\end{bmatrix},\quad x = \begin{bmatrix}f_-\\g_-\\f_+\\g_+\\c_+\\d_+\end{bmatrix},\quad\lambda_k=\begin{bmatrix}\mathcal{R}\left[\bar{\lambda}_{-,k}\right]\\\mathcal{I}\left[\bar{\lambda}_{-,k}\right]\\\mathcal{R}\left[\bar{\lambda}_{+,k}\right]\\\mathcal{I}\left[\bar{\lambda}_{+,k}\right]\end{bmatrix},\;k=0,1.
\]
% The vector $\lambda_k$ is similarly obtained by stacking the real and imaginary parts of $\bar{\lambda}_{-,k}$ and $\bar{\lambda}_{+,k}$, $k=0,1$.
Then
\begin{align*}
\Theta_0 &= \begin{bmatrix}-M_-\\\bm{0}\end{bmatrix},\quad
\Theta_1 = \begin{bmatrix}- M_{\textrm{f}}\\\bm{0}\end{bmatrix}\\
X_0 &=
\begin{bmatrix}\bm{I}&\bm{0}&\bm{0}\\
\bm{0}&\bm{I}&-M_+
\end{bmatrix},\quad
X_1 =
\begin{bmatrix}\bm{I}&\bm{0}&\bm{0}\\
\bm{0}&\bm{I}&-M_{\textrm{f}}
\end{bmatrix}\\
H_0 &= H_1=\bm{I},
\end{align*}
where $\bm{0}$ and $\bm{I}$ are appropriately dimensioned zero and identity matrices. Lastly, we set $Q=\bm{I}$.

% We can now parameterize (\ref{sysnucon}). With a slight abuse of notation, arrange the variables as
% \begin{align*}
% \theta &= \begin{bmatrix}\mu\\\nu\end{bmatrix},\quad x = \begin{bmatrix}f_-\\g_-\\f_+\\g_+\\c_+\\d_+\end{bmatrix}.
% \end{align*}
% Then
% \begin{align*}
% \Theta_0 &= \begin{bmatrix}-M_-\\\bm{0}_2\end{bmatrix},\quad
% \Theta_1 = \begin{bmatrix}- M_{\textrm{f}}\\\bm{0}_2\end{bmatrix}\\
% X_0 &=
% \begin{bmatrix}\bm{I}_2&\bm{0}_2&\bm{0}_2\\
% \bm{0}_2&\bm{I}_2&-M_+
% \end{bmatrix},\quad
% X_1 =
% \begin{bmatrix}\bm{I}_2&\bm{0}_2&\bm{0}_2\\
% \bm{0}_2&\bm{I}_2&-M_{\textrm{f}}
% \end{bmatrix}\\
% H_0 &= H_1=\bm{I}_4.
% \end{align*}
% Lastly, we set $Q=\bm{I}_2$.

% \subsection{Concave restriction}
% We obtain a convex restriction by replacing the bilinear constraint (\ref{bilinear}) with the convex quadratic constraint
% \begin{align}
% 1\leq\sum_{k=0}^1\frac{1}{2}\begin{bmatrix}
% \beta_k\\\theta
% \end{bmatrix}^{\top}\left(\Psi_k-\psi_k\bm{I}\right)\begin{bmatrix}
% \beta_k\\\theta
% \end{bmatrix} - \delta_k.
% \end{align}
% This seems to make the problem infeasible.

%\ifCLASSOPTIONcaptionsoff
%  \newpage
%\fi

\bibliographystyle{IEEEtran}
\bibliography{MainBib,JATBib}
\end{document}